\documentclass[11pt,psfig]{article}
\usepackage{graphicx}
\usepackage{subfigure}

\begin{document}
\bibliographystyle{plain}

\begin{center}

{\Large\bf Uniformly Generating Distribution Functions \\ for Discrete
Random Variables}

\vspace{0.3cm}

{\em Bruno Caprile}

\vspace{0.3cm}

ITC-irst -- Centro per la Ricerca Scientifica e Tecnologica

I-38050 Povo, Trento

Italy

\end{center}

\vspace{1.0cm}

\begin{center}
{\bf Abstract}
\end{center}

An algorithm is presented which, with optimal efficiency, solves the
problem of uniform random generation of distribution functions for an
$n$-valued random variable.

\vspace{1.0cm}

\section{{Introduction} \label{sec:introduction}}

In the general framework of Probabilistic Inference the case occurs
that either for experimental or empirical validation purposes one
needs to generate unbiased collections of distribution functions for
some discrete random variable \cite{LiDam94}. In this note, an
algorithm is presented which efficiently solves the problem.

\section{{The Problem} \label{sec:the-problem}}

Let $x$ be a discrete random variable whose outcomes belong to a
finite set of elementary events $\Omega$, and let $n$ indicate the
cardinality of $\Omega$. Let ${\cal I}_{\Omega}$ be the totality of
distribution functions for $x$. 

\vspace{0.3cm}

\noindent
{\bf Problem:} {\em find an algorithm to sample ${\cal I}_{\Omega}$
uniformly and independently of $n$.}

\vspace{0.3cm}

\noindent
In the following, we shall rely on the existence of a subroutine,
${\cal A}$, able to return series of pseudo random numbers uniformly
distributed in the interval $[0, 1]$. Existence of such subroutine is
thoroughly discussed in \cite{Knuth81}.

\vspace{0.3cm}

\noindent
Let us start by observing that ${\cal I}_{\Omega}$ is naturally
parametrized by $n$ numbers, $x_1, \dots, x_n$, satisfying the
conditions:
\begin{equation}
0 \leq {x_i} \leq 1, ~~~\forall i \in \{1, \dots, n\}
\label{eq-condition-1}
,
\end{equation}

\noindent
and
\begin{equation}
\sum_{i = 1}^{n} x_i = 1
\label{eq-condition-2}
.
\end{equation}

${\cal I}_{\Omega}$ is therefore the $( n - 1 )$-simplex, $S^{n - 1}$,
and our problem in equivalent to finding an algorithm for the uniform
sampling of $S^{n - 1}$. It may be worth reminding that the $n$-volume
of the $n$-simplex tends to zero (super)-exponentially in $n$. This
implies that the na{\"\i}ve sampling strategy consisting in generating
points within the unit $n$-cube, and discarding those falling outside
$S^{n}$ is virtually inapplicable -- even for very small values of
$n$. Other approaches such as that of generating points within the
unit $n$-cube, and rescale them as to satisfy condition $\sum_{i =
1}^{n} x_i = 1$ are plainly wrong.


\section{{Solution} \label{sec:solution}}

Here is the basic idea: for each sample point to be generated on the
$(n - 1)$-simplex, and for each of its first $n - 1$ coordinates,
$x_1, \dots, x_{n - 1}$, randomly sample interval $[0, 1]$ according
to a density function able -- in average -- to assign to each $x_i$
``just its fair share'' of the total amount $\sum_{i = 1}^{n} x_i =
1$. It does not take much to get convinced that such density function
indeed exist for any component $x_{j}$: it is the marginal
distribution of $x_{j}$ over the simplex $S^{n-1}$, given the outcomes
of $x_{1}, \dots, x_{j - 1}$.

\vspace{0.2cm}

The proposed algorithm therefore runs as it follows:

\begin{description}

	\item[1] set $r_{1} = 1$;

	\item[2] set $j = 1$;
	
	\item[3] until $j = n - 1$ 

	\begin{description}

		\item[3.1] randomly extract $x_{j}$ from $[0, r_{j}]$
according to the marginal distribution of $x_{j}$ over the simplex
$S^{n-1}$, given outcomes $x_{1} = \bar{x}_{1}, \dots, x_{j - 1} =
\bar{x}_{j - 1}$, that is according to:

\begin{equation}
\psi(x) = Prob~(x_{j} = x~|~ x_{1} = \bar{x}_{1}, x_{2} =
\bar{x}_{2}, \dots, x_{j -1} = \bar{x}_{j - 1})
\label{eq:psi}
;
\end{equation}

		\item[3.2] set $r_{j + 1} = r_{j} - x_{j}$;

		\item[3.3] set $j = j + 1$;

	\end{description}	

	\item[4] set $\bar{x}_{n} = r_{n}$;

	\item[5] output $(\bar{x}_{1}, \dots, \bar{x}_{n})$.

\end{description} 

Step 3.1 is the crucial one. To perform it, we need to: (1) determine
$\psi$ for any $n$ and any set of outcomes, $\bar{x}_{1}, \dots,
\bar{x}_{j - 1}$; (2) sample interval $[0, r_{j}]$ according to
$\psi$.

\subsection{{Determining $\psi(x)$} \label{subsec:determining-psi(x)}}

\begin{figure}
\includegraphics[scale=.4]{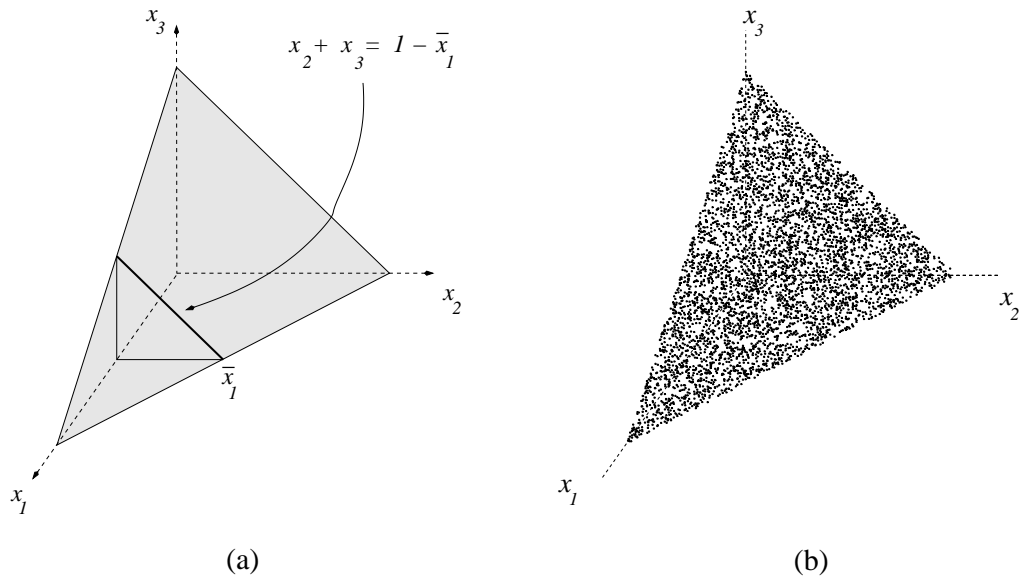}
\caption{(a) {\em The 3-dimensional case: for any outcome of $x_{1}$
{\em rescaling} of the 1-simplex is determined. The
marginal distribution of $x_{1}$ is therefore proportional to the
1-volume of $S^{1}$.} (b) {\em 5000 samples of $S^{2}$ as
obtained by applying the proposed algorithm.}}
\label{fig:compose-marginal-with-5000-samples}
\end{figure}

Let us start by observing that the marginal distribution of $x_{1}$
must be proportional to the $(n - 2)$-volume of the subset of $R^{n}$
defined by $x_{2} + x_{3} + \dots + x_{n} = 1 - x_{1}$. Let us
indicate such subset with $S_{x_{1}}$. For any $x_{1}$, $S_{x_{1}}$ is
just a {\em rescaling} of the $(n - 2)$-simplex, and its volume is
therefore proportional to $(1 - x_{1})^{n - 2}$ (see
Fig. \ref{fig:compose-marginal-with-5000-samples}a). The marginal
distribution of $x_{1}$ can therefore be written in the form
$\psi(x_{1}) = \alpha (1 - x_{1})^{n - 2}$, where factor $\alpha$ is
determined via the normalization condition:

\begin{equation}
\alpha \int_{0}^{1} (1 - x_{1})^{n - 2} dx_{1} = 1
\label{eq:psi-normalization-x1}
.
\end{equation}

\noindent
This yields $\alpha = n - 1$, and then:

\begin{equation}
\psi(x_{1}) = ( n - 1 ) (1 - x_{1})^{n - 2} 
\label{eq:psi-x1}
.
\end{equation}

\vspace{0.3cm}

The same process can be iterated for all the other components $x_{j}$,
$1 < j < n$, accounting for the fact that any $x_{j}$ is now to be
limited to the range $[0, r_{j}]$. This implies that $\psi(x_{j})$
must be proportional to $(r_{j} - x_{j})^{n - 2}$, and it is an
interesting fact that dependence of $\psi(x_{j})$ from outcomes
$\bar{x}_{1}, \dots, \bar{x}_{j - 1}$ is just contained in their sum
$1 - r_{j}$. Thus, we can write:

\begin{equation}
\beta \int_{0}^{r_{j}} (r_{j} - x_{j})^{n - 2} dx_{j} = 1
\label{eq:psi-normalization}
,
\end{equation}

\noindent
which yields $\beta = \frac{n - 1}{r_{j}^{n - 1}}$, and, finally

\begin{equation}
\psi(x_{j}) = \frac{n - 1}{r_{j}^{n - 1}} (1 - x_{j})^{n - 2}
\label{eq:psi-xj}
.
\end{equation}

\vspace{0.3cm}

\noindent
The cumulative function of the marginal distribution of $x_{j}$ is
then:

\begin{equation}
\Psi(x_{j}) = 1 - \left(\frac{r_{j} - x_{j}}{r_{j}}\right)^{n - 1}
\label{eq:cumulative}
.
\end{equation}

\subsection{{Sampling $[0, s]$ according to $\psi(x)$}
\label{subsec:sampling-[0,s]-according-to-psi(x)}} 

As it is well know \cite{PreTeuVetFla92}, sampling a random variable
$x$ according to a given distribution function, $\psi(x)$, is readily
obtained once that the inverse of the cumulative function of $x$,
$\Psi^{-1}$, is known. Indeed, Eq. \ref{eq:cumulative}, guarantees
that, for any $n$ and $j$:

\begin{equation}
\Psi^{-1}(\xi) = r_{j} [1 - (1 - \xi)^{\frac{1}{n - 1}}]
\label{eq:inverse-cumulative}
.
\end{equation}

\section{{Efficiency} \label{sec:Efficiency}}

The proposed algorithm is optimally efficient: in dimension $n$ it
requires just $n - 1$ runs of subroutine ${\cal A}$, plus $n - 2$
calls of function $\Psi^{-1}$, whose complexity is constant in $n$.

\paragraph{Acknowledgements}

The author wish to thank B. Walsh for signaling an error in formula
(9) as reported in an earlier version of the preprint.

\end{document}